\documentclass[12pt,oneside,english]{amsart}
\usepackage[T1]{fontenc}
\usepackage[latin1]{inputenc}
\usepackage{geometry}
\usepackage{amssymb}

\makeatletter
 \theoremstyle{plain}    
 \newtheorem{thm}{Theorem}[section]
 \numberwithin{equation}{section} 
 \numberwithin{figure}{section} 
 \theoremstyle{plain}
 \theoremstyle{plain}    
 \newtheorem{prop}[thm]{Proposition} 
 \theoremstyle{plain}    
 \newtheorem{lem}[thm]{Lemma} 
 \theoremstyle{remark}
 \newtheorem{rem}[thm]{Remark}

\def\makebbb#1{
    \expandafter\gdef\csname#1\endcsname{
        \ensuremath{\Bbb{#1}}}
}
\makebbb{R}
\makebbb{N}
\makebbb{Z}
\makebbb{C}
\makebbb{G}
\makebbb{E}
\makebbb{H}
\makebbb{P}
\makebbb{K}
\makebbb{B}

\def\hpq0{h^{p,q}_{\leq 0}}
\def\Hpq0{\H_{\leq 0}^{p,q}}

\def\dbar{\bar\partial}
\def\ddbar{\partial\dbar}

\def\C{{\mathbb C}}

\def\H{{\mathcal H}}

\def\Re{{\rm Re\,  }}

\def\be{\begin{equation}}
\def\ee{\end{equation}}

\usepackage{babel}
\makeatother
\begin{document}

\title{A direct approach to Bergman kernel asymptotics for positive line
bundles }

\email{robertb@math.chalmers.se bob@math.chalmers.se \newline \quad\quad\quad
johannes@math.polytechnique.fr }

\author{Robert Berman, Bo Berndtsson, Johannes Sjöstrand}

\begin{abstract}
We give an elementary proof of the existence of an asymptotic expansion
in powers of $k$ of the Bergman kernel associated to $L^{k}$, where
$L$ is a positive line bundle over a compact complex manifold. We
also give an algorithm for computing the coefficients in the expansion.
\end{abstract}
\maketitle

\section{Introduction\label{sec:Introduction}}

Let $L$ be a holomorphic line bundle with a positively curved hermitian
metric $\phi$, over a complex manifold $X.$ Then $i/2$ times the
curvature form $\partial\overline{\partial}\phi$ of $L$ defines
a Kähler metric on $X,$ that induces a scalar product on the space
of global sections with values in $L.$ The orthogonal projection
$P$ from $L^{2}(X,L)$ onto $H^{0}(X,L),$ the subspace of holomorphic
sections, is the Bergman projection. Its kernel with respect to the
scalar product is the Bergman kernel $K$ of $H^{0}(X,L)$; it is
a section of $\overline{L}\boxtimes L$ over $X\times X.$ It can
also be characterized as a reproducing kernel for the Hilbert space
$H^{0}(X,L),$ i.e \begin{equation}
\alpha(x)=(\alpha,K_{x})\label{eq:repr prop intro}\end{equation}
 for any element $\alpha$ of $H^{0}(X,L),$ where $K_{x}=K(\cdot,x)$
is identified with a holomorphic section of $L\otimes\overline{L_{x}},$
where $L_{x}$ denotes the fiber of $L$ over $x$. The restriction
of $K$ to the diagonal is a section of $\overline{L}\otimes L$ and
we let $B(x)=\left|K(x,x)\right|$ be its pointwise norm.

Even though the Bergman kernel is of course impossible to compute
in general, quite precise asymptotic formulas, when we replace $L$
by $L^{k}$ and $\phi$ by $k\phi$, are known, see Zelditch \cite{Zelditch}
and Catlin \cite{Catlin}. Namely: \be
K_x(y)e^{-k\psi(y,x)}=\frac{k^n}{\pi^n}( 1+ \frac{b_1(x,y)}{k}+\frac{b_2(x,y)}{k^2}
+...) \ee
as $k$ goes to infinity.

Here $\psi$ is an (almost) holomorphic extension of $\phi$ and the
$b_{j}$s are certain hermitian functions. In particular, on the diagonal
$y=x$ we have an asymptotic series expansion for \[
K_{x}(x)e^{-k\phi(x)}=\frac{k^{n}}{\pi^{n}}(1+\frac{b_{1}(x,x)}{k}+\frac{b_{2}(x,x)}{k^{2}}+...).\]
 The functions $b_{j}(x,x)$ contain interesting geometric information
of $X$ with the Kahler metric $i\ddbar\phi/2$, see Lu \cite{Lu}.

In \cite{Zelditch} and \cite{Catlin} the existence of an asymptotic
expansion is proved using a formula, due to Boutet de Monvel and Sjöstrand,
for the boundary behaviour of the Bergman-Szegö kernel for a strictly
pseudoconvex domain, \cite{Boutet-Sjöstrand}, extending an earlier
result of C Fefferman, \cite{Fefferman}, to include also the off-diagonal
behaviour. The purpose of this paper is to give a direct proof of
the existence of this expansion, the main point being that it is actually
simpler to construct an asymptotic formula directly. Our method also
gives an effective way of computing the terms $b_{j}$ in the asymptotic
expanison. Even though the inspiration for the construction comes
from the calculus of Fourier integral operators with complex phase,
the arguments in this paper are elementary.

The method of proof uses localization near an arbitrary point of $X$.
Local holomorphic sections of $L^{k}$ in a small coordinate neighbourhood
, $U$, are just holomorphic functions on $U$, and the local norm
is a weighted $L^{2}$-norm over $U$ with weight function $e^{-k\phi}$
where $\phi$ is a strictly plurisubharmonic function. Using the ideas
from \cite{Sjöstrand} we then compute \textit{local asymptotic Bergman
kernels} on $U$. These are holomorphic kernel functions, and the
scalar product with such a kernel function reproduces the values of
holomorphic functions on $U$ up to an error that is small as $k$
tends to infinity. Assuming the bundle is globally positive it is
then quite easy to see that the global Bergman kernel must be asymptotically
equal to the local kernels.

Many essential ideas of our approach were already contained in the
book \cite{Sjöstrand} written by the third authour. Here we use them
in order to find a short derivation of the Bergman kernel asymptotics.
For the closely related problem of finding the Bergman kernel for
exponentially weighted spaces of holomorphic functions, this was done
by A. Melin and the third author \cite{Melin-Sjöstrand}, but in the
present paper we replace a square root procedure used in that paper
by a more direct procedure, which we think is more convenient for
the actual computations of the coefficients in the asymptotic expansions.
There are also close relations to the subject of weighted integral
formulas in complex analysis \cite{Berndtsson-Andersson}. We have
tried to make the presentation almost self-contained, hoping that
it may serve as an elementary introduction to certain micro-local
techniques with applications to complex analysis and differential
geometry.

\section{\label{sec:The-local-asymptotic}The local asymptotic Bergman kernel}

The local situation is as follows. Fix a point in $X.$ We may choose
local holomorhic coordinates $x$ centered at the point and a holomorphic
trivialization $s$ of $L$ such that \begin{equation}
\left|s\right|^{2}=e^{-\phi(x)},\label{eq:def of local fiber
  m}\end{equation}
 where $\phi$ is a smooth real valued function. $L$ is positive
if and only if all local functions $\phi$ arising this way are strictly
plurisubharmonic. We will call $\phi_{0}(x)=\left|x\right|^{2}$ the
\emph{model fiber metric}, since it may be identified with the fiber
metric of a line bundle of constant curvature on $\C^{n}.$ The Kähler
form, $\omega$, of the metric on our base manifold $X$ is given
by $i/2$ times the curvature form of $L$, \[
\omega=i\partial\bar{\partial}\phi/2.\]
 The induced volume form on $X$ is equal to $\omega_{n}:=\omega^{n}/n!$.
Now any local holomorphic section $u$ of $L^{k}$ may be written
as $us^{\otimes k}$ where $u$ is a holomorphic function. The local
expression of the norm of a section of $L^{k}$ over $U$ is then
given by \[
\left\Vert u\right\Vert _{k\phi}^{2}:=\int_{U}\left|u\right|^{2}e^{-k\phi}\omega_{n}\]
 where $u$ is a holomorphic function un $U.$ The class of all such
functions $u$ with finite norm is denoted by $H_{k\phi}(U).$

We now turn to the construction of local asymptotic Bergman kernels.
The main idea is that since a posteriori Bergman kernels will be quite
concentrated near the diagonal, we require a local Bergman kernel
to satisfy the reproducing formula (1.1) locally, up to an error which
is exponentially small in $k$. In the sequel we fix our coordinate
neighbourhood to be the unit ball of $\mathbb{C}^{n}$. Let $\chi$
be a smooth function supported in the unit ball $B$ and equal to
one on the ball of radius $1/2$. We will say that $K_{k}$ is a \emph{reproducing
kernel} mod $O(e^{-\delta k})$ for $H_{k\phi}$ if for any fixed
$x$ in some neighbourhood of the origin we have that for any local
holomorphic function $u_{k},$\begin{equation}
u_{k}(x)=(\chi u_{k},K_{k,x})_{k\phi}+O(e^{k(\phi(x)/2-\delta)})\left\Vert u\right\Vert _{k\phi},\label{eq:local repre prop}\end{equation}
 uniformly in some neighbourhood of the origin. Furthermore, if $K_{k,x}(y)$
is holomorphic in $y$ we say that $K_{k,x}$ is a \emph{Bergman kernel}
mod $O(e^{-\delta k})$.

Given a positive integer $N,$ Bergman and reproducing kernels mod
$O(k^{-N})$ are similarly defined.

\subsection{\label{sub:The-model-case}Local reproducing kernels mod $e^{-\delta k}$ }

Let $\phi$ be a strictly plurisubharmonic function in the unit ball.
and let $u$ be a holomorphic function in the ball such that \[
\Vert u\Vert^{2}:=\int_{B}|u|^{2}e^{-k\phi}<\infty.\]
 We shall first show that (cf \cite{Sjöstrand}) integrals of the
form \begin{equation}
c_{n}(k/2\pi)^{-n}\int_{\Lambda}e^{k\theta\cdot(x-y)}u(y)d\theta\wedge dy,\label{int}\end{equation}
 define reproducing kernels mod $e^{-\delta k}$ for suitably choosen
\textit{contours} \[
\Lambda=\{(y,\theta);\theta=\theta(x,y)\}.\]
 Here we think of $x$ as being fixed (close to the origin) and let
$y$ range over the unit ball, so that $\Lambda$ is a $2n$-dimensional
submanifold of $B_{y}\times\C_{\theta}^{n}$, and $c_{n}=i^{n}(-1)^{n(n+1)/2}=i^{-n^{2}}$
is a constant of modulus 1 chosen so that $c_{n}d\bar{y}\wedge dy$
is a positive form. Let us say that such a contour is \textit{good}
if uniformly on $\Lambda$ for $x$ in some neighbourhood of the origin
and $|y|\leq1$ \[
2\Re\theta\cdot(x-y)\leq-\delta|x-y|^{2}-\phi(y)+\phi(x).\]
 Note that, by Taylor's formula \[
\phi(x)=\phi(y)+2\Re\sum Q_{i}(x,y)(x_{i}-y_{i})+\sum\phi_{i\bar{j}}(x_{i}-y_{i})\overline{(x_{j}-y_{j})}+o(x-y)^{2},\]
 where $\sum Q_{i}(x,y)(x_{i}-y_{i})$ is the part of the second order
Taylor expansion which is holomorphic in $x$. Hence, if $\phi$ is
strictly plurisubharmonic, \[
\theta(x,y)=Q(x,y)\]
 is a good contour, depending on $x$ in a holomorphic way. In particular,
$\theta=\bar{y}$ defines a good contour for $\phi(x)=|x|^{2}$. 

\begin{prop}
For any good contour, \[
u(x)=(k/2\pi)^{n}c_{n}\int_{\Lambda}e^{k\theta\cdot(x-y)}u(y)\chi(y)d\theta\wedge dy+O(e^{k(\phi(x)/2-\delta)})\Vert u\Vert_{k\phi},\]
 for $x$ in some fixed neighbourhood of 0 if $u$ is an element of
$H_{k\phi}(B)$. 
\end{prop}
\begin{proof}
For $s$ a real variable between 0 and $\infty$, we let \[
\Lambda_{s}=\{(x,y,\theta);\theta+s(\bar{x}-\bar{y})\in\Lambda\},\]
 and denote by $\eta=\eta_{k}$ the differential form \[
\eta=c_{n}(k/2\pi)^{n}e^{k\theta\cdot(x-y)}u(y)\chi(y)d\theta\wedge dy.\]
 Our presumptive reproducing formula is the integral, $I_{0}$, of
$\eta$ over $\Lambda_{0}$ and it is easy to see that the limit of
\[
I_{s}:=\int_{\Lambda_{s}}\eta\]
 as $s$ goes to infinity equals $u(x)$. (This is because $c_{n}(s/2\pi)^{n}e^{-s|x-y|^{2}}d\bar{y}\wedge dy$
tends to a Dirac measure at $x$ as $s$ tends to infinity.) The difference
between $I_{0}$ and $I_{s}$ is by Stokes formula \[
I_{0}-I_{s}=\int_{B\times[0,s]}dh^{*}(\eta)\]
 where $h$ is the homotopy map \[
h(y,\lambda)=(y,\theta(x,y)-\lambda(\bar{x}-\bar{y})).\]
 Now, \[
d\eta=c_{n}(k/2\pi)^{n}e^{k\theta\cdot(x-y)}ud\chi\wedge d\theta\wedge dy.\]
 This equals 0 if $|y|<1/2$, and since $\theta$ is good we have
the estimate \[
|dh^{*}(\eta)|\leq Ck^{n}e^{k(-(\delta/2+\lambda)|x-y|^{2}-\phi(y)/2+\phi(x)/2)}(1+\lambda)^{n}|u(y)|.\]
 If $|x|$ is, say smaller than 1/4, $|x-y|\geq1/4$ when $d\eta$
is different from 0, so we get \[
|\int dh^{*}(\eta)|\leq Ck^{n}e^{k(\phi(x)/2-\delta)}\int_{|y|>1/2}|u(y)|e^{-k\phi(y)/2}\int_{0}^{s}(1+\lambda)^{n}e^{-k\lambda}d\lambda\]
 with a smaller $\delta$. By the Cauchy inequality the first integral
in the right hand side is dominated by \[
\Vert u\Vert_{k\phi}.\]
 Since the last integral is bounded by a constant independent of $k$
we get the desired estimate. 
\end{proof}
Thus we have a family of reproducing kernels mod $e^{-\delta k}$.
When $\phi=|y|^{2}$ and $\theta=\bar{y}$, the kernel in the representation
- $e^{\bar{x}\cdot y}$ - is also holomorphic in $y$ so we even have
an asymptotic Bergman kernel mod $e^{-\delta k}$. To achieve the
same thing for general weights we need to introduce a bit more flexibility
in the construction by allowing a more general class of amplitudes
in the integral. We will replace the function $\chi$ in the integral
by $\chi(1+a)$ for a suitably choosen function $a$, where $a$ has
to be chosen to give an exponentially small contribution to the integral.

For this we consider differential forms \[
A=A(x,y,\theta,k)=\sum A_{j}(x,y,\theta,k)\widehat{d\theta_{j}}\]
 of bidegree $(n-1,0)$. By $\widehat{d\theta_{j}}$ we mean the wedge
product of all the differentials $d\theta_{i}$ except $d\theta_{j}$,
with a sign chosen so that $d\theta_{j}\wedge\widehat{d\theta_{j}}=d\theta$.
We assume that $A$ has an asymptotic expansion of order 0 \[
A\sim A_{0}+k^{-1}A_{1}+...\]
 By this we mean that for any $N\geq0$ \[
A-\sum_{0}^{N}A_{m}k^{-m}=O(k^{-N-1})\]
 uniformly as $k$ goes to infinity.

We assume also that the coefficients are holomorphic (in the smooth
case almost holomorphic) for $x$, $y$ and $\theta$ of norm smaller
than 2. Let \[
ad\theta=e^{-k\theta\cdot(x-y)}d_{\theta}e^{k\theta\cdot(x-y)}A,\]
 so that \begin{equation}
a=D_{\theta}\cdot A+k(x-y)\cdot A=:\nabla A,\end{equation}
 where $D_{\theta}=\partial/\partial\theta$. We will say that a function
$a$ arising in this way is a negligible amplitude. In the applications
we will also need to consider finite order approximations to amplitude
functions. Let \[
A^{(N)}=\sum_{0}^{N}A_{m}/k^{m}\]
 and similarily \[
a^{(N)}=\sum_{0}^{N}a_{m}/k^{m}.\]
 Then \[
a^{(N)}=\nabla A^{(N+1)}-D_{\theta}\cdot A_{N+1}/k^{N+1},\]
 so $a^{(N)}$ is a negligible amplitude modulo an error term which
is $O(1/k^{N+1})$.

\begin{prop}
For any good contour $\Lambda$ and any negligible amplitude $a$,
\[
u(x)=c_{n}(k/2\pi)^{n}\int_{\Lambda}e^{k\theta\cdot(x-y)}u(y)\chi(y)(1+a)d\theta\wedge dy+O(e^{k(\phi(x)/2-\delta)})\Vert u\Vert_{k\phi},\]
 for all $x$ in a sufficiently small neighbourhood of the origin
if $u$ is an element of $H_{k\phi}(B)$. Moreover \[
u(x)=c_{n}(k/2\pi)^{n}\int_{\Lambda}e^{k\theta\cdot(x-y)}u(y)\chi(y)(1+a^{(N)})d\theta\wedge dy+O(e^{k\phi(x)/2}/k^{N+1-n})\Vert u\Vert_{k\phi},\]

\end{prop}
\begin{proof}
For the first statement we need to verify that the contribution from
$a$ is exponentially small as $k$ tends to infinity. But \[
\int_{\Lambda}e^{k\theta\cdot(x-y)}u(y)\chi(y)ad\theta\wedge dy=\int_{\Lambda}u(y)\chi(y)d_{\theta}(e^{k\theta\cdot(x-y)}A)\wedge dy=\]
 \[
=\int_{\Lambda}\chi d\left(u(y)e^{k\theta\cdot(x-y)}A\wedge dy\right)=-\int_{\Lambda}d\chi\wedge u(y)e^{k\theta\cdot(x-y)}A\wedge dy.\]
 Again, the last integrand vanishes for $|y|<1/2$ and is, since $\Lambda$
is good, dominated by a constant times \[
|u(y)|e^{k(-\delta|x-y|^{2}-\phi(y)/2+\phi(x)/2)}\]
 The last integral is therefore smaller than \[
\Vert u\Vert O(e^{k(\phi(x)/2-\delta)})\]
 so the first formula is proved. The second formula follows since
by the remark immediately preceeding the proposition, $a^{(N)}$ is
a good amplitude modulo an error of order $1/k^{N+1}$. 
\end{proof}
The condition that an amplitude function $a$ can be written in the
form (2.6) can be given in an equivalent very useful way. For this
we will use the infinite order differential operator \[
Sa=\sum_{0}^{\infty}\frac{1}{(k)^{m}(m!)}(D_{\theta}\cdot D_{y})^{m}.\]
 This is basically the classical operator that appears in the theory
of pseudodifferential operators when we want to replace an amplitude
$a(x,y,\theta)$ by an amplitude $b(x,\theta)$ independent of $y$,
see \cite{Hörmander1}. We let $S$ act on $(n-1)$-forms $A$ as
above componentwise. We say that $Sa=b$ for $a$ and $b$ admitting
asymptotic expansions if all the coefficients of the powers $(1/k)^{m}$
in the expansion obtained by applying $S$ to $a$ formally equal
the corresponding coefficients in the expansion of $b$. No convergence
of any kind is implied. That $Sa$ equals $b$ to order $N$ means
that the same thing holds for $m\leq N$. Note also that since formally
\[
S=e^{D_{\theta}\cdot D_{y}/k},\]
 we have that \[
S^{-1}=e^{-D_{\theta}\cdot D_{y}/k}=\sum_{0}^{\infty}\frac{1}{(-k)^{m}(m!)}(D_{\theta}\cdot D_{y})^{m}.\]

\begin{lem}
Let \[
a\sim\sum(1/k)^{m}a_{m}(x,y,\theta)\]
 be given. Then there exists an $A$ satisfying (2.6) asymptotically
if and only if \[
Sa|_{x=y}=0.\]
 Moreover the last equation holds to order $N$ if and only if $a^{(N)}$
can be written \be
a^{(N)}=\nabla A^{(N+1)} +O(1/k^{N+1}). \ee

\end{lem}
\begin{proof}
Note first that $S$ commutes with $D_{\theta}$ and that \[
S\left((x-y)\cdot A\right)=(x-y)\cdot SA-(1/k)D_{\theta}\cdot SA.\]
 Moreover \[
\nabla A=D_{\theta}\cdot A+k(x-y)\cdot A,\]
 so it follows that \be
S\nabla A= SD_\theta\cdot A +kS(x-y)\cdot A= \ee
\[
D_{\theta}SA+k(x-y)\cdot SA-D_{\theta}S\cdot A=k(x-y)\cdot SA.\]
 Thus, if $a$ admits a representation $a=\nabla A$, then $Sa$ must
vanish for $x=y$. Similarily, if \[
a^{(N)}=\nabla A^{(N+1)}+O(1/k^{N+1}).\]
 it follows that $Sa^{(N)}|_{y=x}=0$ to order $N$.

Conversely, assume $Sa|_{y=x}=0$ . Then $Sa=(x-y)\cdot B$ for some
form $B$. But (2.6) implies that \[
\nabla S^{-1}=kS^{-1}(x-y)\cdot\]
 so \[
a=S^{-1}\left((x-y)\cdot B\right)=(1/k)\nabla S^{-1}B\]
 and (2.4) holds with $A=1/kS^{-1}B$. If the equation $Sa|_{y=x}=0$
only holds to order $N$, then \[
Sa^{(N)}=(x-y)\cdot B^{(N)}\]
 to order $N$. Hence \[
a^{(N)}=S^{-1}\left((x-y)\cdot B^{(N)}\right)=(1/k)\nabla S^{-1}B^{(N)}=(1/k)\nabla(S^{-1}B)^{(N)}\]
 to order $N$, so (2.7) holds with $A^{(N+1)}=1/k(S^{-1}B)^{(N)}$.
\end{proof}

\subsection{The phase}

Let us now see how to choose the contour $\Lambda$ to get the phase
function $\psi(x,\bar{y})$ appear in the expression \[
e^{\theta\cdot(x-y)}.\]
 In this section we still assume that the plurisubharmonic function
$\phi$ is real analytic and let $\psi(x,y)$ be the unique holomorphic
function of $2n$ variables such that \[
\psi(x,\bar{x})=\phi(x).\]
 By looking at the Taylor expansions of $\psi$ and $\phi$ one can
verify that \begin{equation}
2\Re\psi(x,\bar{y})-\phi(x)-\phi(y)\leq-\delta|x-y|^{2}\end{equation}
 for $x$ and $y$ sufficiently small. Following an idea of Kuranishi,
see \cite{Hörmander2}, \cite{Grigis-Sjöstrand}, we now find a holomorphic
function of $3n$ variables, $\theta(x,y,z)$ that solves the division
problem \begin{equation}
\theta\cdot(x-y)=\psi(x,z)-\psi(y,z),\label{eq:imposing change variable}\end{equation}
 This can be done in many ways, but any choice of $\theta$ satisfies
\[
\theta(x,x,z)=\psi_{x}(x,z).\]
 To fix ideas, we take \[
\theta(x,y,z)=\int_{0}^{1}\partial\psi(tx+(1-t)y,z)dt\]
 with $\partial$ denoting the differential of $\psi$ with respect
to the first $n$ variables.

Since $\theta(x,x,z)=\psi_{x}(x,z)$ it follows that \[
\theta_{z}(0,0,0)=\psi_{xz}(0,0)=\phi_{x\bar{x}}(0,0)\]
 is a nonsingular matrix. Therefore \[
(x,y,z)\rightarrow(x,y,\theta)\]
 defines a biholomorphic change of coordinates near the origin. After
rescaling we may assume that $\psi$ is defined and satisfies (2.7)
and that the above change of coordinates is well defined when $|x|$,
$|y|$ and $|z|$ are all smaller than 2. We now define $\Lambda$
by \[
\Lambda=\{(y,\theta);z=\bar{y}\}.\]
 Thus, on $\Lambda$, $\theta$ is a holomorphic function of $x$,
$y$ and $\bar{y}$. The point of this choice is that by (2.8), on
$\Lambda$, \[
\theta\cdot(x-y)=\psi(x,\bar{y})-\psi(y,\bar{y}).\]
 Therefore we get the right phase function in our kernel and by (2.7)
\[
2\Re\theta\cdot(x-y)=2\Re\psi(x,\bar{y})-2\phi(y)\leq\phi(x)-\phi(y)-\delta|x-y|^{2},\]
 which means that $\Lambda$ is a good contour in the sense of the
previous section. By Proposition 2.2 we therefore get the following
proposition, where we use the notation $\beta$ for the standard Kähler
form in $\mathbb{C}^{n}$, \[
\beta=i/2\sum dy_{j}\wedge d\bar{y}_{j}.\]

\begin{prop}
\label{pro:integral correct phase}Suppose that $u$ is in $H_{k\phi}.$
If $a(x,y,\theta,1/k)$ is a negligible amplitude, we have \begin{equation}
u(x)=\end{equation}
 \[
=(k/\pi)^{n}\int\chi_{x}e^{k(\psi(x,\overline{y})-\psi(y,\overline{y}))}(\det\theta_{\bar{y}})u(y)(1+a)\beta_{n}+\]
 \[
+O(e^{k(\phi(x)/2+\delta)})\left\Vert u\right\Vert _{k\phi},\]
 with $a=a(x,y,\theta(x,y,\bar{y})$. Moreover \[
u(x)=\]
 \[
=(k/\pi)^{n}\int\chi_{x}e^{k(\psi(x,\overline{y})-\psi(y,\overline{y}))}(\det\theta_{\bar{y}})u(y)(1+a^{(N)})\beta_{n}+\]
 \[
+O(e^{k\phi(x)/2}k^{n-N-1})\left\Vert u\right\Vert _{k\phi}\]

\end{prop}

\subsection{The amplitude}

In order to get an asymptotic Bergman kernel from (2.10) we need to
choose the amplitude $a$ so that \[
\det\theta_{\bar{y}}(1+a)=B(x,\bar{y})\det\psi_{y\bar{y}},\]
 with $B$ analytic. Polarizing in the $y$-variable, i e replacing
$\bar{y}$ by $z$, this means that \[
(1+a(x,y,\theta(x,y,z,1/k))=B(x,z,1/k)\det\psi_{yz}(y,z)/\det\theta_{z}(x,y,z),\]
 where $B$ is a analytic and independent of $y$. Consider this as
an equation between functions of the variables $x$, $y$ and $\theta$.
Let \[
\Delta_{0}(x,y,\theta)=\det\psi_{yz}(y,z)/\det\theta_{z}(x,y,z)=\det\partial_{\theta}\psi_{y}.\]
 Since $\psi_{y}=\theta$ when $y=x$ we have that $\Delta_{0}=1$
for $y=x$. We need $a$ to be representable in the form (2.6) which
by the previous lemma means that $Sa=0$ for $y=x$. Equivalently,
$S(1+a)=1$ for $y=x$, so we must solve \be
S\left(B(x,z(x,y,\theta),1/k)\Delta_0(x,y,\theta)\right) =1
\ee
for $y=x$. This equation should hold in the sense of formal power
series which means that the coefficient of $1/k^{0}$ must equal 1,
whereas the coefficient of each power $1/k^{m}$ must vanish for $m>0$.
In the computations $x$ is held fixed and $z=z(y,\theta)$. The first
equation is \be
b_0(x,z(x,y,\theta))\Delta_0(x,x,\theta) =1. \ee
 This means that $b_{0}(x,z(x,\theta))=1$, for all $\theta$ which
implies that $b_{0}$ is identically equal to 1.

The second condition is \be
(D_\theta\cdot D_y)\left(b_0\Delta_0\right)+b_1\Delta_0=0
\ee
for $y=x$. Since we already know that $b_{0}=1$ this means that
\[
b_{1}(z(x,\theta))=-(D_{\theta}\cdot D_{y})\left(\Delta_{0}\right)|_{y=x},\]
 which again determines $b_{1}$ uniquely. Continuing in this way,
using the recursive formula

\be
\sum_0^m\frac{(D_\theta\cdot D_y)^l}{l!} \left(b_{m-l}\Delta_0\right)|_{x=y}=0
\ee
for $m>0$ we can determine all the coefficients $b_{m}$, and hence
$a$. Then $Sa|_{y=x}=0$ so $Sa^{(N)}|_{y=x}=0$ to order $N$, and
the next proposition follows from Propositions 2.4 and 2.3.

\begin{prop}
Suppose that $\phi$ is analytic. Then there are analytic functions
$b_{m}(x,z)$ defined in a fixed neighbourhood of $x$ so that for
each $N$\begin{equation}
(k/\pi)^{n}(1+b_{1}(x,\overline{y})k^{-1}+...+b_{N}(x,\overline{y})k^{-N})e^{k\psi(x,\overline{y})},\label{eq:exp in prop}\end{equation}
 is an asymptotic Bergman kernel mod $O(k^{-(N+1)})$. 
\end{prop}

\subsection{Computing $b_{1}$}

Let us first recall how to express some Riemannian curvature notions
in Hermitian geometry. The Hermitian metric two-form $\omega:=\frac{i}{2}H_{ij}dy^{i}\wedge\overline{dy^{j}}$
determines a connection $\eta$ on the complex tangent bundle $TX$
with connection matrix (with respect to a holomorphic frame) \begin{equation}
\eta=H^{-1}\partial H=:\sum\eta_{j}dy_{j}.\label{eq:connection}\end{equation}
 The curvature is the matrix valued two-form $\overline{\partial}\eta$
and the scalar curvature $s$ is $\Lambda\textrm{$\textrm{Tr}$}\overline{\partial}\eta$
where $\Lambda$ is contraction with the metric form $\omega.$ Hence,
in coordinates centered at $x$ where $H(0)=I$ the scalar curvature
$s$ at $0$ is given by \begin{equation}
s(0)=-\textrm{Tr}(\sum\frac{\partial}{\partial\overline{y_{j}}}\eta_{j}),\label{eq:s}\end{equation}
 considering $\eta$ as matrix. We now turn to the computation of
the coefficient $b_{1}$ in the expansion \ref{eq:exp in prop}. By
the definition of $\theta$ we have that\begin{equation}
\theta_{i}(x,y,z)=\psi_{y_{i}}(y,z)+\frac{1}{2}\sum_{k}(\frac{\partial}{\partial y_{k}}\psi_{y_{i}})(y,z)(x^{k}-y^{k})+...\label{eq:exp of theta}\end{equation}
 Differentiating with respect to $z$ gives\[
\theta_{z}=H+\frac{1}{2}\partial_{y}H(x-y)+....,\]
 where $H=H(y,z).$ Multiplying both sides by $H^{-1}$ and inverting
the relation we get \begin{equation}
\theta_{z}^{-1}H=I-\frac{1}{2}(H^{-1}\partial H)(x-y)+....,\label{eq:quotient}\end{equation}
 Taking the determinant of both sides in formula \ref{eq:quotient}
gives\begin{equation}
\Delta_{0}=1-\frac{1}{2}\textrm{Tr}\eta(x-y)+....,\label{eq:m is connec}\end{equation}
 Hence, equation (2.14) now gives, since $-\frac{\partial}{\partial y}(x-y)=1,$
that\[
b_{1}(0,0)=(\frac{\partial}{\partial\theta}\cdot(-\frac{1}{2}\textrm{Tr}\eta)_{x=y}=-\frac{1}{2}\frac{\partial}{\partial\overline{y}}\cdot\textrm{Tr}\eta\]
 showing that $b_{1}(x,\overline{x})=\frac{s}{2},$ according to \ref{eq:s}.

\subsection{Twisting with a vector bundle $E$}

We here indicate how to extend the previous calculation to the case
of sections with values in $L^{k}\otimes E,$ where $E$ is a holomorphic
vector bundle with a hermitian metric $G$ (see also \cite{Keller}).
First observe that $u(x)$ is now, locally, a holomorphic vector and
the Bergman kernel may be identified with a matrix $K(x,\overline{y})$
such that \[
u(x)=\int K(x,\overline{y})G(y,\overline{y})u(y)\psi_{y\overline{y}}e^{-k\psi(y,\overline{y}))}d\bar{y}\wedge dy.\]
 To determine $K$ one now uses the ansatz \[
K(x,\overline{y})=c_{n}(k/2\pi)^{n}e^{k(\psi(x,\overline{y})}B(x,\overline{y},k^{-1})G(x,\overline{y})^{-1}.\]
 Then the condition on the amplitude function becomes \begin{equation}
(1+a(x,y,\theta(x,y,z),1/k)\det(\frac{\partial\theta}{\partial z}(x,y,z))=\end{equation}
 \[
=B(x,z,1/k)G(x,z)^{-1}G(y,z)\det(\psi_{yz}).\]
 where $a$ now is a matrix valued form , i.e. $\Delta_{0}$ in section
2.3 is replaced by the matrix $\Delta_{G}:=\Delta_{0}G(x,z)^{-1}G(y,z).$
Note that \[
G(x,z)^{-1}G(y,z)=I-G^{-1}(y,z)\frac{\partial}{\partial y}G(y,z)(x-y)+...=:I-\eta_{E}(y,z)(x-y)+...,\]
 where $\eta_{E}:=G^{-1}\frac{\partial}{\partial y}G$ is the connection
matrix of $E.$ Hence, the equation \ref{eq:m is connec} is replaced
by \[
\Delta_{G}=1\textrm{-(Tr}\eta/2\otimes I+\eta_{E})(x-y)+...\]
 The same calculation as before then shows that the matrix $b_{1}(0,0)$
is given by\[
b_{1}(0,0)=-\frac{1}{2}\frac{\partial}{\partial\overline{y}}Tr\eta\otimes I-\frac{\partial}{\partial\overline{y}}\cdot\eta_{E}=\frac{s}{2}\otimes I+\Lambda\Theta_{E},\]
 where $\Theta_{E}:=\overline{\partial}\eta_{E}$ is the curvature
matrix of $E$ and $\Lambda$ denotes contraction with the metric
two-form $\omega.$

\textbf{Remark:} Let $K_{k}$ be the Bergman kernel of $H^{0}(X,L^{k}),$
defined with respect a general volume form $\mu_{n}.$ Then the function
$G:=\mu_{n}/\omega_{n}$ defines a hermitian metric on the trivial
line bundle $E$ and the asymptotics of $K_{k}$ can then be obtained
as above.

\subsection{\label{sub:l smooth}Smooth metrics.}

Denote by $\psi(y,z)$ any almost holomorphic extension of $\phi$
from $\overline{\Delta}=\{ z=\bar{y}\}$, i.e. an extension such that
the anti-holomorphic derivatives vanish to infinite order on $\overline{\Delta}.$
We may also assume that $\overline{\psi(y,\overline{z})}=\psi(z,\overline{y}).$
That $\psi$ is almost holomorphic means that for any multi index
$\alpha$ \begin{equation}
D^{\alpha}(\overline{\partial}\psi)=0,\label{eq:prop of psi almost}\end{equation}
(where $D^{\alpha}$ is the local real derivative of order $\alpha)$
when evaluated at a point in $\overline{\Delta},$ i.e. when $z=\bar{y}.$
Let now \begin{equation}
\theta=\int_{0}^{1}(\partial_{y}\psi)(tx+(1-t)y,z)dt,\,\theta^{*}=\int_{0}^{1}(\overline{\partial}_{y}\psi)(tx+(1-t)y,z)dt,\label{eq:def of smooth change of v}\end{equation}
(where $\partial_{y}\psi$ denotes the vector of partial holomorphic
derivatives w.r.t the first argument of $\psi)$ so that \begin{equation}
(x-y)\theta+\overline{(x-y)}\theta^{*}=\psi(x,z)-\psi(y,z)\label{eq:not kuranishi}\end{equation}
 Then the smooth map corresponding to $(x,y,z)\mapsto(x,y,\theta)$
is locally smoothly invertible for the same reason as in the analytic
case, since $\overline{\partial}_{z}\theta=0$ when $x=y=\overline{z}.$
Define the algebra $\mathcal{A}$ of all functions almost holomorphic
when $x=y=\overline{z}$ as the set of smooth functions, $f$, of
$x$ $y$ and $z$, such that \[
D^{\alpha}\dbar f=0\]
 for all multi indices $\alpha,$ when $x=y=\bar{z}$. For a vector
valued function we will say that it is in $\mathcal{A}$, if its components
are in $\mathcal{A}.$ We also define the vanishing ideal $\mathcal{I}^{\infty}$
as the set of smooth functions $f$ such that \[
D^{\alpha}f=0,\]
 for all multi indices $\alpha,$ when $x=y=\bar{z}$. Hence, if $f$
belongs to $\mathcal{A}$ then (the coefficients of) $\dbar f$ will
belong to $\mathcal{I}^{\infty}$.

Note that $\psi(tx+(1-t)y,z)$ is in $\mathcal{A}$ for each fixed
$t.$ Hence $\theta$ is in $\mathcal{A}$ and $\theta^{*}$ is in
$\mathcal{I}^{\infty},$ so that \ref{eq:not kuranishi} gives \begin{equation}
(x-y)\theta=\psi(x,\overline{y})-\psi(y,\overline{y})+O(\left|x-y\right|^{\infty})\label{eq:exponent in smooth}\end{equation}
when $\bar{z}=y.$ The next simple lemma is used to show that the
contribution of elements in $\mathcal{I}^{\infty}$ to the phase function
and the amplitude is negligable.

\begin{lem}
\label{lem:negligable}Let  $f_{i}$ be elements of the vanishing ideal
$\mathcal{I}^{\infty}$ 
and let $b(x,y)$ be a local smooth function.
Then\[
\int\chi_{x}(y)e^{k(\psi(x,\overline{y})-\psi(y,\overline{y})+f_{1}(x,y,\overline{y}))}(b(x,y)+f_{2}(x,y,\overline{y}))u(y)\beta_{n}(y)=\]
\[
=\int\chi_{x}(y)e^{k(\psi(x,\overline{y})-\psi(y,\overline{y})}b(x,y)u(y)\beta_{n}(y)+O(k^{-\infty})\left\Vert u\right\Vert _{k\phi}\]

\end{lem}
\begin{proof}
First observe that if $f_{i}$ is an element of $\mathcal{I}^{\infty},$
then $f_{i}(x,y,\overline{y})=O(|x-y|^{\infty}.$ Moreover, we have
\begin{equation}
\left|\left|x-y\right|^{2N}e^{k(\psi(x,\overline{y})-\phi(x)/2-\phi(y)/2)}\right|\leq Ck^{-2N}(k\left|x-y\right|)^{2})^{N}e^{-k\delta\left|x-y\right|^{2}}=O(k^{-2N}),\label{eq:bd}\end{equation}
 where we have used (2.7). Combining this bound with the Cauchy-Schwartz
inequality proves the lemma when $f_{1}=0.$ Now write\[
e^{k(\psi(x,\overline{y})-\psi(y,\overline{y})+f_{1}(x,y,\overline{y}))}=e^{k(\psi(x,\overline{y})-\psi(y,\overline{y})}+\int_{0}^{1}\partial_{t}(e^{(\psi(x,\overline{y})-\psi(y,\overline{y})+tf_{1}(x,y,\overline{y}))})dt.\]
By \ref{eq:bd} the second term gives a contribution which is of the
order $O(k^{-\infty}).$ Hence the general case follows.
\end{proof}
. 

\begin{prop}
\label{thm:local bergmank smooth case}Suppose that $L$ is smooth.
Then there exists an asymptotic reproducing kernel \emph{$K_{k}^{(N)}$
mod $O(k^{n-N-1})$} for $H_{k\phi}$ , such that \be
K_{k}^{(N)}(x,\overline{y})=e^{k\psi(x,\overline{y})}(b_{0}+b_{1}k^{-1}+...+b_{N}k^{-N})\label{form} \ee
 where $b_{i}$ is a polynomial in the derivatives $\partial_{x}^{\alpha}\overline{\partial}_{y}^{\beta}\psi(x,\overline{y})$
of the almost holomorphic extension $\psi$ of $\phi.$ In particular,
\begin{equation}
e^{-k(\phi(x)/2+\phi(y)/2)}(D_{x,y}^{\alpha}(\overline{\partial}_{x},\partial_{y}))K(x,y)=O(k^{-\infty})\label{eq:k almost holom}\end{equation}
 uniformly in $x$ and $y$ for any given $\alpha.$
\end{prop}
\begin{proof}
We go through the steps in the proof of the analytic case and indicate
the necessary modifications.

First we determine the coefficients $b_{m}(x,z)$ in the same way
as in the analytic case, i e by fixing $x$ and solving \[
S(B(z)\Delta_{0})|_{y=x}=1\]
 Here $S$ has the same meaning as before and in particular contains
only derivatives with respect to $\theta$ and no derivatives with
respect to $\bar{\theta}$. The difference is that $\Delta_{0}$ is
no longer analytic so $B$ will not be holomorphic, but it will still
belong to $\mathcal{A}$ since $\Delta_{0}$ does.

We next need to consider lemma 2.3 with $a\in\mathcal{A}.$ Then we
get that \[
a\in\mathcal{A},\,(Sa)_{y=x}=O(k^{-N-1})\Leftrightarrow\exists A\in\mathcal{A}:\, a=\nabla A+O(k^{-N-1})\,\textrm{mod}\mathcal{\, I^{\infty}}\]
 Indeed, this follows from the argument in the analytic case and the
fact that if $c(=c(x,y,z))\in\mathcal{A},$ then \[
c(x,x,z)=0\Leftrightarrow\exists d\in\mathcal{A}:\, c=(x-y)d\,\textrm{mod}\mathcal{\, I^{\infty}}\]
 as can be seen by defining $d$ by \[
d=\int_{0}^{1}(\partial_{y}c)(x,x+(1-t)y,z)dt.\]
Here $\nabla$ also has the same meaning as before and contains only
a derivative with respect to $\theta$ and no derivative with respect
to $\bar{\theta}$. Then Proposition 2.2 holds as before except that
there will be one extra contribution in the application of Stokes
theorem coming from $\dbar_{\theta}A$ (when $z=\bar{y}).$ Since
$\dbar_{\theta}A$ vanishes to infinite order when $x=y=\bar{z}$,
it gives a contribution to the integral which is $O(1/k^{N})$ for
any $N$ by lemma \ref{lem:negligable}.

We therefore get from Proposition 2.2 a reproducing kernel of the
form claimed in \eqref{form} except that the phase function equals
\[
k\theta\cdot(x-y)=k(\psi(x,\bar{y})-\phi(y))+f)\]
 with $f$ in $\mathcal{A}$. Again by lemma \ref{lem:negligable}
we may remove $f$ at the expense of adding a contribution to the
to the integral which is negligable, i.e which is $O(1/k^{N})$ for
any $N$. 
\end{proof}

\section{The global Bergman kernel}

In this section we will show that, if the curvature of $L$ is positive
everywhere on $X,$ then the global Bergman kernel $K_{k}$ of $H^{0}(X,L^{k})$
is asymptotically equal to the local Bergman kernel $K_{k}^{(N)}$
of $H_{k\phi}$ (constructed in section \ref{sec:The-local-asymptotic}).

Recall (section \ref{sec:Introduction}) that the Bergman kernel $K$
associated to $L$ is a section of $\overline{L}\boxtimes L$ over
$X\times X.$ By restriction $K_{x}$ is identified with a holomorphic
section of $\overline{L_{x}}\otimes L,$ where $L_{x}$ is the fiber
of $L$ over $x.$ Given any two vector spaces $E$ and $F,$ the
scalar product on $L$ extends uniquely to a pairing \begin{equation}
(\cdot,\cdot):\, L\otimes E\times L\otimes F\rightarrow E\otimes F,\label{eq:pairing}\end{equation}
 linear over $E$ and anti-linear over $F.$ In terms of this pairing
$K_{y}$ has the global reproducing property \begin{equation}
\alpha(y)=(\alpha,K_{y})\label{eq:repre prop in global}\end{equation}
 for any element $\alpha$ of $H^{0}(X,L).$ By taking $\alpha=K_{x}$
(so that $E=L_{x}$ and $F=\overline{L_{y}}$ in \ref{eq:pairing})
one gets \begin{equation}
K(y,x):=K_{x}(y)=(K_{x},K_{y}).\label{eq:k as scalar
  product}\end{equation}
 This also implies that $\overline{K(x,y)}=K(y,x)$ and that \begin{equation}
K(x,x)=(K_{x},K_{x})=\left\Vert K_{x}\right\Vert ^{2}.\label{eq:b as norm}\end{equation}
 $K(x,x)$ is a section of $\bar{L}\otimes L$. Its norm as a section
to this bundle is the Bergman function, which in a local frame with
respect to which the metric on $L$ is given by $e^{-\phi}$ equals
\[
B(x)=K(x,x)e^{-\phi(x)}.\]
 Notice also that by the Cauchy inequality we have an extremal characterization
of the Bergman function: \[
B(x)=\sup|s(x)|^{2}\]
 where the supremum is taken over all holomorphic sections to $L$
of norm not greater than 1.

We now denote by $K_{k}$ the Bergman kernel associated to $L^{k}$,
and write $B_{k}$ for the associated Bergman function. It follows
from the extremal characterization of the Bergman function and the
submeanvalue inequality for a holomorphic section $s$ over a small
ball with radius roughly $1/k^{1/2}$ that \[
B_{k}\leq Ck^{n},\]
 uniformly on $X$ (see e g \cite{Berman}).

Let now $K_{x}^{(N)}(y)$ be the local Bergman kernel of propositions
2.5-2.6, where the coefficients $b_{m}$ are given by (2.15), \be
\overline{K^{(N)}_x(y)}= (k/\pi)^n(1+b_{1}(x,\overline{y})k^{-1}+...+b_{N}(x,\overline{y})
k^{-N})e^{k\psi(x,\overline{y})}. \ee
By construction, the coefficients $b_{m}(x,z)$ are holomorphic if
the metric on $L$ - locally represented by $\phi$ - is real analytic.
In case $\phi$ is only smooth the $b_{m}$s are almost holomorphic,
meaning that \[
\dbar_{xz}b_{m}\]
 vanishes to infinite order when $z=\bar{x}$.

Replacing $K_{y}$ in the relation \ref{eq:k as scalar product} with
the local Bergman kernel $K_{k}^{(N)}$ will now show that $K_{k}=K_{k}^{(N)}$
up to a small error term.

\begin{thm}
Assume that the smooth line bundle $L$ is globally positive. Let
$K_{k}^{(N)}$ be defined by (3.5), where the coefficients $b_{m}$
are determined by the recursion (2.15).

If the distance $d(x,y)$ is sufficiently small, then \begin{equation}
K_{k}(x,y)=K_{k}^{(N)}(x,y)+O(k^{n-N-1})e^{k(\phi(x)/2+\phi(y)/2)}\textrm{,}\label{eq:statement  
  of th 
  global}\end{equation}
 Moreover, \[
D^{\alpha}(K_{k}(x,y)-K_{k}^{(N)}(x,y))=O(k^{m+n-N-1})e^{k(\phi(x)/2+\phi(y)/2)}\]
 if $D^{\alpha}$ is any differential operator with respect to $x$
and $y$ of order at most $m$. 
\end{thm}
\begin{proof}
Let us first show that \begin{equation}
K_{k}(y,x)=(\chi K_{k,x},K_{k,y}^{(N)})+O(k^{n-N-1})e^{k(\phi(x)/2+\phi(y)/2)}\label{eq:pr glob berg claim a}\end{equation}
 where $\chi$ is a cut-off function equal to 1 in a neighbourhood
of $x$ which is large enough to contain $y$. Fixing $x$ and applying
Proposition 2.5 to $u_{k}=K_{k,x}$ gives \ref{eq:pr glob berg claim a}
with the error term \[
e^{\phi(y)/2}O(k^{-N-1})\left\Vert K_{x}\right\Vert .\]
 Now, by \ref{eq:b as norm} and the estimate for $B_{k}$ \[
\left\Vert K_{k,x}\right\Vert ^{2}=B_{k}(x)e^{k\phi(x)}\leq Ck^{n}e^{k\phi(x)},\]
 This proves \ref{eq:pr glob berg claim a} with uniform convergence.

Next we estimate the difference \[
u_{k,y}(x):=K_{k,y}^{(N)}(x)-(\chi K_{y}^{(N)},K_{k,x}).\]
 Since the scalar product in this expression is the Bergman projection,
\[
P_{k}(\chi K_{k,y}^{(N)})(x),\]
 $u_{k,y}$ is the $L^{2}$-minimal solution to the $\dbar$-equation
\[
\dbar u_{k,y}=\dbar(\chi K_{k,y}^{(N)}).\]
 The right hand side equals \[
(\dbar\chi)K_{k,y}^{(N)}+\chi\dbar K_{k,y}^{(N)}.\]
 Since $\chi$ equals 1 near $y$ it follows from (2.9) and the explicit
form of $K_{k,y}^{(N)}$ that the first term is dominated by \[
e^{-\delta k}e^{k(\phi(\cdot)/2+\phi(y)/2)}.\]
 The second term vanishes identically in the analytic case. In the
smooth case $\dbar K_{k}^{(N)}$ can by proposition \ref{thm:local bergmank smooth case}
be estimated by \be
O(1/k^\infty)e^{k(\phi(\cdot)/2+\phi(y)/2)}. \ee
Altogether $\dbar u_{k,y}$ is therefore bounded by (3.7), so by the
Hörmander $L^{2}$-estimate we get that \[
\Vert u_{k,y}\Vert^{2}\leq O(1/k^{\infty})e^{k\phi(y)/2}.\]
 But, since the estimate on $\dbar u_{k,y}$ is even uniform, we get
by a standard argument involving the Cauchy integral formula in a
ball around $x$ of radius roughly $1/k^{1/2}$ that $u_{k,y}$ satisfies
a pointwise estimate \[
|u_{k,y}(x)|^{2}\leq O(1/k^{\infty})e^{k(\phi(y)/2+\phi(x)/2)}.\]

Combining this estimate for $u_{k,y}(x)$ with (3.6) we finally get
\be
|K^{(N)}_{k,y}(x)-\overline{K_k(y,x)}|e^{-k\phi(x)/2-k\phi(y)/2}\leq
O(1/k^{N+1}). \ee
Since $K_{k}$ is hermitian ( i e $K_{k}(x,y)=\overline{K_{k}(y,x)}$)
this proves the proposition except for the statement on convergence
of derivatives.

In the analytic case the convergence of derivatives is, by the Cauchy
estimates, an automatic consequence of the uniform convergence, since
the kernels are holomorphic in $x$ and $\bar{y}$. In the smooth
case, we have that \[
\dbar K_{k}^{(N)}(x,\bar{z})=O(1/k^{\infty})e^{k(\phi(\cdot)/2+\phi(y)/2)}.\]
 This implies that the Cauchy estimates still hold for the difference
between $K_{k}$ and $K_{k}^{(N)}$, up to an error which is $O(1/k^{\infty})$,
and so we get the convergence of derivatives even in the smooth case.
\end{proof}
\begin{rem}
The proof above actually shows that the asymptotic expansion for the
global Bergman kernel $K_{k}(x,y)$ holds close to any point $x$
which is in $X(0)$ (the open subset of $X$ where the curvature form
of $\phi$ is positive) and is such that $x$ satisfies the following
global condition: for any given $\overline{\partial}-$closed $(0,1)-$form
$g_{k}$ with values in $L^{k}$ supported in some fixed neighbourhood
of $x$ we may find sections $u_{k}$ with values in $L^{k}$ such
that\begin{equation}
\overline{\partial}u_{k}=g_{k},\label{eq:inhom dbar}\end{equation}
 on $X$ and \begin{equation}
\left\Vert u_{k}\right\Vert _{k\phi}\leq C\left\Vert g_{k}\right\Vert _{k\phi}\label{eq:horm est2}\end{equation}
After this paper was written this observation was used in \cite{Berman2}
to obtain an asymptotitc expansion of $K_{k}(x,y)$ on a certain subset
of $X(0)$ for any Hermitian line bundle $(L,\phi)$ over a projective
manifold $X.$
\end{rem}

\end{document}